
\input graphicx.tex
\input amssym.def
\input amssym
\magnification=1000
\baselineskip = 0.20truein
\lineskiplimit = 0.01truein
\lineskip = 0.01truein
\vsize = 8.8truein
\voffset = 0.1truein
\parskip = 0.09truein
\parindent = 0.3truein
\settabs 12 \columns
\hsize = 6.0truein
\hoffset = 0.21truein

\setbox\strutbox=\hbox{%
\vrule height .708\baselineskip
depth .292\baselineskip
width 0pt}
\font\caps=cmcsc10
\font\bigtenrm=cmr10 at 14pt

\def\sqr#1#2{{\vcenter{\vbox{\hrule height.#2pt
\hbox{\vrule width.#2pt height#1pt \kern#1pt
\vrule width.#2pt}
\hrule height.#2pt}}}}
\def\square{\mathchoice\sqr46\sqr46\sqr{3.1}6\sqr{2.3}4}

\centerline{\bigtenrm FINDING DISJOINT SURFACES IN 3-MANIFOLDS}
\tenrm
\vskip 14pt
\centerline{MARC LACKENBY}
\vskip 18pt
\centerline{\caps 1. Introduction}
\vskip 6pt

In 3-manifold theory, properly embedded surfaces play a key role. 
It is particularly interesting when a compact 3-manifold $M$ contains
{\sl two} properly embedded disjoint surfaces $S_1$ and $S_2$ such that
$M - (S_1 \cup S_2)$ is connected. It is this scenario
that we will investigate in this paper. Our aim is to show that 
the 3-manifolds where this situation occurs are actually very plentiful.
Moreover, we will see how algebraic methods can be profitably
used to detect the existence of such a pair of surfaces. In particular,
the infinite dihedral group ${\Bbb Z}_2 \ast {\Bbb Z}_2$
will play a central role. This is because a compact connected 3-manifold $M$
contains a pair of such surfaces if and only if $\pi_1(M)$ admits 
a surjective homomorphism onto ${\Bbb Z}_2 \ast {\Bbb Z}_2$
(see Theorem 3.1.)

Note that we do not require that $S_1$ and $S_2$ are 2-sided.
The existence of two disjoint properly embedded 2-sided surfaces $S_1$ and $S_2$
in a 3-manifold $M$ with $M - (S_1 \cup S_2)$ connected appears to be a much more rare occurrence, 
and is equivalent to the existence of a surjective homomorphism $\pi_1(M) \rightarrow {\Bbb Z} \ast {\Bbb Z}$.
(See [1] for example, where obstructions to the existence of such a homomorphism are given.)

The following is our main result.

\noindent {\bf Theorem 1.1.} {\sl Let $M$ be a compact connected orientable
3-manifold, and suppose that $\partial M$ is non-empty and contains no 2-spheres. Then
$M$ contains two properly embedded disjoint surfaces $S_1$ and $S_2$
such that $M - (S_1 \cup S_2)$ is connected if and only if
$M$ is neither a ${\Bbb Z}_2$ homology solid torus nor a ${\Bbb Z}_2$
homology cobordism between two tori.}

\noindent {\bf Corollary 1.2.} {\sl Let $L$ be a link in $S^3$. Then the exterior $X$ of $L$ contains
two properly embedded disjoint surfaces $S_1$ and $S_2$ such that $X - (S_1 \cup S_2)$
is connected if and only if one of the following holds:
\item{(i)} $L$ has at least three components, or
\item{(ii)} $L$ has two components, which have even linking number.

}

It is an interesting and not completely straightforward exercise
to construct such surfaces in the exterior of the
Whitehead link. We will do so explicitly in Section 4.

The case of 2-component links was studied by Hillman in [2], where
he stated Corollary 1.2 in this case. (It appears in the middle of
the second full paragraph on page 176 in [2].)

The plan of the paper is as follows. In Section 2, we prove the foundational
result that the existence of a surjective homomorphism from a finitely generated
group $G$ onto ${\Bbb Z}_2 \ast {\Bbb Z}_2$ is equivalent to the existence of an index
2 subgroup $K$ of $G$ with $b_1(K) > b_1(G)$. In Section 3, we first prove
that, for a compact connected 3-manifold $M$, the existence of a surjective homomorphism
$\pi_1(M) \rightarrow {\Bbb Z}_2 \ast {\Bbb Z}_2$ is equivalent to the
existence of two properly embedded disjoint surfaces as in the statement 
of Theorem 1.1. We then use these two facts to prove one direction of
Theorem 1.1. The starting point is the well-known inequality
$b_1(M) \geq b_1(\partial M)/2$ for any compact orientable 3-manifold $M$.
Assuming that $M$ is a 3-manifold as in Theorem 1.1 which is neither 
a ${\Bbb Z}_2$ homology solid torus nor a ${\Bbb Z}_2$
homology cobordism between two tori,
we find a double cover $\tilde M$ of $M$ such that $b_1(\partial \tilde M) > b_1(\partial M)$,
and, after some work, we deduce that $b_1(\tilde M) > b_1(M)$.
The results in Section 2 and earlier in Section 3 then give the required surfaces.
The existence of these two surfaces is thereby proved, but the construction
is group-theoretic and far from explicit. In Section 4, we give a more
geometric way of finding these surfaces, which can be used easily in practice.
We examine the case of the exterior of the Whitehead link, which is
an instructive example. In Section 5,
we introduce methods from profinite group theory. We show that the
existence of a surjective homomorphism from a finitely generated group
$G$ to ${\Bbb Z}_2 \ast {\Bbb Z}_2$ is detected by the pro-$2$ completion of $G$.
These profinite techniques are used to prove the other direction of Theorem 1.1,
which establishes that certain 3-manifolds as described in the theorem do not
contain two disjoint properly embedded surfaces whose union is non-separating.
In Section 6, we also use profinite group theory to control the ${\Bbb Z}_2$ homology classes
of the surfaces for certain 3-manifolds. For example, we show that in the case
where $M$ is a compact orientable 3-manifold that has the same ${\Bbb Z}_2$
homology as a handlebody (other than a solid torus), then $S_1$ and $S_2$ may
be chosen to represent any pair of distinct non-trivial classes in $H_2(M, \partial M; {\Bbb Z}_2)$. 
The key input here is the fact that, in this case, $\pi_1(M)$ has the same
pro-2 completion as a non-abelian free group. In Section 7, we pose some
questions, which may stimulate further research in this area.

\vskip 18pt
\centerline{\caps 2. Surjections to the infinite dihedral group}
\vskip 6pt

We start with the following group-theoretic result.

\noindent {\bf Theorem 2.1.} {\sl Let $G$ be a finitely generated group,
and let $K$ be an index $2$ subgroup. Then the following are equivalent:
\item{(i)} There is a surjective homomorphism $\phi \colon G \rightarrow {\Bbb Z}_2 \ast {\Bbb Z}_2$
such that $K$ is the kernel of $\pi \phi$, where $\pi \colon {\Bbb Z}_2 \ast {\Bbb Z}_2 \rightarrow
{\Bbb Z}_2$ is the homomorphism which is an isomorphism on each factor.
\item{(ii)} $b_1(K) > b_1(G)$.

}

Note that this has the following immediate corollary.

\noindent {\bf Corollary 2.2.} {\sl Let $G$ be a finitely generated group. Then the following are
equivalent:
\item{(i)} There is a surjective homomorphism $G \rightarrow {\Bbb Z}_2 \ast {\Bbb Z}_2$.
\item{(ii)} For some index two subgroup $K$ of $G$, $b_1(K) > b_1(G)$.

}

\noindent {\sl Proof of Theorem 2.1.} (i) $\Rightarrow$ (ii). 
Suppose that there is a surjective homomorphism $\phi \colon G \rightarrow {\Bbb Z}_2 \ast {\Bbb Z}_2$
such that $K$ is the kernel of $\pi \phi$. We claim that $b_1(K) > b_1(G)$. 

Now, the inclusion $i \colon K \rightarrow G$ induces a 
homomorphism $i^\ast \colon H^1(G; {\Bbb R}) \rightarrow H^1(K; {\Bbb R})$.
This is clearly an injection. This is because $H^1(G; {\Bbb R})$ may be viewed as the set of homomorphisms $G \rightarrow {\Bbb R}$
and if a homomorphism $G \rightarrow {\Bbb R}$
is zero when restricted to $K$, then it is zero on all of $G$.

We will show that in fact this injection $i^\ast$ is not a surjection. This will prove that
$b_1(K) > b_1(G)$ as required. 

The kernel of the homomorphism $\pi$ is a subgroup $A$ of 
${\Bbb Z}_2 \ast {\Bbb Z}_2$ that is infinite cyclic.
Let $\psi$ be the composition
$K \rightarrow A \rightarrow {\Bbb R}$, where the first map is the
restriction of $\phi$ to $K$ and the second map is the standard inclusion
of the infinite cyclic group into ${\Bbb R}$. We will show that $\psi$
is not in the image of $i^\ast$.

For suppose that $\psi$ is in this image. This means that $\psi$ extends to
a homomorphism $\tilde \psi \colon G \rightarrow {\Bbb R}$. Let $a$ and $b$
be generators of the factors of ${\Bbb Z}_2 \ast {\Bbb Z}_2$.  Let $g_a$ and
$g_b$ be elements of $G$ that are sent by $\phi$ to these generators. Then $g_a g_b$ lies in
$K$ and its image under $\psi$ is $1$, say. On the other hand, $g_b g_a$ also
lies in $K$ and its image under $\psi$ is $-1$, because $ab$ and $ba$ 
are inverses. But, because ${\Bbb R}$ is an abelian group,
$\psi(g_a g_b) = \tilde \psi(g_a g_b) = \tilde \psi(g_b g_a) = \psi(g_b g_a)$,
which is a contradiction.

(ii) $\Rightarrow$ (i). Let $K$ be an index 2 subgroup of $G$ with
$b_1(K) > b_1(G)$.
Let $V$ be the vector space $H^1(K; {\Bbb R})$. This acted
on by $G/K = {\Bbb Z}_2$ via the conjugation action of $G$ on $K$. Let
$\tau$ be the automorphism of $H^1(K; {\Bbb R})$ induced by
the non-trivial element of $G/K$.
Since $\tau$ is an involution, it gives rise to a decomposition $V = V_- \oplus V_+$,
where $V_-$ and $V_+$ are the $-1$ and $+1$ eigenspaces of $\tau$. 

We claim that the 
eigenspace $V_+$ consists of precisely the homomorphisms $K \rightarrow {\Bbb R}$
that extend to $G$. For suppose that a homomorphism $\psi \colon K \rightarrow {\Bbb R}$
extends to a homomorphism $\tilde \psi \colon G \rightarrow {\Bbb R}$. Let $g$ be an element
of $G - K$. The image of  $\psi$ under $\tau$ sends $k \in K$ to $\psi(g^{-1} k g) = 
\tilde \psi (g^{-1} k g) = \tilde \psi(k) = \psi(k)$. So, $\psi$ lies in $V_+$. Conversely,
suppose that $\psi \colon K \rightarrow {\Bbb R}$ lies in $V_+$ and is therefore invariant
under $\tau$. Then we extend $\psi$ to $\tilde \psi \colon G \rightarrow {\Bbb R}$ by
defining
$$\eqalign{
\tilde \psi(k) &= \psi(k) \cr
\tilde \psi(gk) &= \psi(k) + (\psi(g^2))/2 \cr
}$$
Here, $g$ is a fixed element of $G - K$, and $k$ is an arbitrary element of $K$. It is easy to 
check that $\tilde\psi$ is a homomorphism, which clearly extends $\psi$. This proves the claim.

So, the dimension of $V_+$ is $b_1(G)$. 
Now, we are assuming that $b_1(K) > b_1(G)$, and so $V_-$ is therefore
non-zero. We can view $H^1(K; {\Bbb Z})$ as a lattice in $H^1(K; {\Bbb R})$,
consisting of those homomorphisms $K \rightarrow {\Bbb R}$ that have image in
${\Bbb Z}$. The action of $G/K$ on $H^1(K; {\Bbb R})$ leaves this
lattice invariant as a set, simply because it sends a homomorphism
$K \rightarrow {\Bbb Z}$ to another such homomorphism. Thus, we may
find non-zero elements of $H^1(K; {\Bbb Z}) \cap V_-$, as follows. Take any element
$\alpha$ of $H^1(K; {\Bbb Z})$ not in $V_+$, and consider $\alpha - \tau \alpha$,
where $\tau \alpha$ is the image of $\alpha$ under the action of $\tau$.
Thus, as we know that $H^1(K; {\Bbb Z}) \cap V_-$ is non-zero, we
may find a primitive element $\psi$ in $H^1(K; {\Bbb Z}) \cap V_-$.
This corresponds to a surjective homomorphism $\psi \colon K \rightarrow {\Bbb Z}$.
We now use this to define a surjective homomorphism $\phi \colon G \rightarrow {\Bbb Z}_2 \ast {\Bbb Z}_2$.
Pick any element $g$ in $G - K$. Note that $g^2 \in K$ and that
$\psi(g^2) = \psi(g^{-1} g^2 g) = - \psi(g^2)$, and hence $\psi(g^2) = 0$.
Define a function
$$\eqalign{
G & \buildrel \phi \over \longrightarrow {\Bbb Z}_2 \ast {\Bbb Z}_2 \cr
k &\mapsto (ab)^{\psi(k)}, \cr
gk &\mapsto a (ab)^{\psi(k)}. \cr}$$
Here $k$ is an arbitrary element of $K$.
This is easily checked to be a homomorphism. For example, suppose that
$k_1$, $k_2 \in K$. We check that $\phi(k_1) \phi(g k_2) = \phi(k_1 g k_2)$:
$$\phi(k_1)\phi(gk_2) = (ab)^{\psi(k_1)} a (ab)^{\psi(k_2)}$$
$$\phi(k_1 g k_2) = \phi( g (g^{-1} k_1 g) k_2) = a (ab)^{-\psi(k_1)} (ab)^{\psi(k_2)}
= a (ba)^{\psi(k_1)} (ab)^{\psi(k_2)}.$$
\vfill \eject
We also check that $\phi(g k_1) \phi(g k_2) = \phi(g k_1 g k_2)$:
$$\phi(gk_1)\phi(gk_2) = a (ab)^{\psi(k_1)} a (ab)^{\psi(k_2)} = (ab)^{-\psi(k_1) + \psi(k_2)}$$
$$\phi(gk_1 g k_2) = \phi(g^2 (g^{-1} k_1 g) k_2) = 
(ab)^{\psi(g^2)} (ab)^{-\psi(k_1)} (ab)^{\psi(k_2)}.$$
This homomorphism $\phi$ is surjective. This can be seen as follows. Since
$\psi$ is surjective, there is some $k \in K$ such that $\psi(k) = 1$.
So, $\phi(gk) = aab = b$, and $\phi(kgk) = \phi(k) \phi(gk) = abb = a$.

Finally, note that ${\rm Ker}(\pi) = A$ consists of precisely those elements
of ${\Bbb Z}_2 \ast {\Bbb Z}_2$ of the form $(ab)^m$ for some $m \in {\Bbb Z}$.
So, the kernel of $\pi \phi$ is exactly $K$, as required. This proves the theorem. $\square$

\vskip 18pt
\centerline{\caps 3. Constructing disjoint surfaces in 3-manifolds}
\vskip 6pt

The following is fairly well known. See, for example [3], where a related result is proved.

\noindent {\bf Theorem 3.1.} {\sl Let $M$ be a compact connected 3-manifold. Then, the
following are equivalent:
\item{(i)} There is a surjective homomorphism $\pi_1(M) \rightarrow {\Bbb Z}_2 \ast {\Bbb Z}_2$.
\item{(ii)} There are two disjoint properly embedded surfaces $S_1$ and $S_2$ in $M$
such that $M - (S_1 \cup S_2)$ is connected.

}

\noindent {\sl Proof.} (ii) $\Rightarrow$ (i). Suppose that there are surfaces $S_1$ and $S_2$ as in (ii).
We will use these to
define a continuous map $f \colon M \rightarrow ({\Bbb R}P^3) \vee ({\Bbb R}P^3)$ such that $f_\ast \colon
\pi_1(M) \rightarrow \pi_1(({\Bbb R}P^3) \vee ({\Bbb R}P^3)) = {\Bbb Z}_2 \ast {\Bbb Z}_2$ is
surjective.

Let $N(S_1)$ and $N(S_2)$ be disjoint regular neighbourhoods of $S_1$
and $S_2$. Then $N(S_i)$ is an $I$-bundle over $S_i$ in which $S_i$ lies as a zero-section. 
We now define a map $f_i \colon N(S_i) \rightarrow {\Bbb R}P^3$. 

Pick a cell structure on $S_i$. This lifts to a cell structure on $\tilde S_i = {\rm cl}(\partial N(S_i) - \partial M)$,
which is the $(\partial I)$-bundle over $S_i$.
This extends to a cell structure on $N(S_i)$, as follows. The fibre over each 0-cell
of $S_i$ becomes a 1-cell of $N(S_i)$. The interior of each 1-cell of $S_i$ has
inverse image in $N(S_i) - \tilde S_i$ that is an open disc, which we declare to
be the interior of a 2-cell. Similarly, each 2-cell of $S_i$ induces a
3-cell of $N(S_i)$.
We define $f_i$ one cell at a time, starting with the 0-cells, then the 1-cells and so on.
Give ${\Bbb R}P^3$ the usual cell structure, with one 0-cell, one 1-cell, one 2-cell and one 3-cell.
Now define $f_i$ on the 0-cells of $N(S_i)$, by sending them to the unique 0-cell of
${\Bbb R}P^3$. Map each 1-cell of $N(S_i)$ that misses $S_i$ also to the 0-cell of 
${\Bbb R}P^3$. Send each 1-cell of $N(S_i)$ that intersects $S_i$ around the 1-cell
of ${\Bbb R}P^3$ so that the interior of the cell is mapped in homeomorphically.
There are two different ways of doing this. We pick one arbitrarily. The 2-cells of $N(S_i)$
come in two varieties. There are 2-cells that lie in $\tilde S_i$. We send these
to the 0-cell of ${\Bbb R}P^3$. The 2-cells that are vertical in $N(S_i)$ we map to the
2-cell of ${\Bbb R}P^3$. This is possible because the boundary of each vertical 2-cell of $N(S_i)$ runs
over two vertical 1-cells, and so its boundary has been mapped to a loop in
${\Bbb R}P^3$ that is homotopically trivial. Finally, the 3-cells of $N(S_i)$ may be
mapped in because $\pi_2({\Bbb R}P^3) = 0$.

Thus, we have defined a map $f_i \colon N(S_i) \rightarrow {\Bbb R}P^3$. Note that
the image of $\tilde S_i$ is the 0-cell of ${\Bbb R}P^3$.

Now form the wedge ${\Bbb R}P^3 \vee {\Bbb R}P^3$ by gluing the two copies of 
${\Bbb R}P^3$ along the 0-cells. We may define $f \colon M \rightarrow {\Bbb R}P^3 \vee {\Bbb R}P^3$ 
by sending points outside of $N(S_1) \cup N(S_2)$ to the basepoint of the wedge,
and by mapping in $N(S_i)$ to the $i$th copy of ${\Bbb R}P^3$ using $f_i$. 

We claim that $f_\ast \colon \pi_1(M) \rightarrow \pi_1(({\Bbb R}P^3) \vee ({\Bbb R}P^3))$ is
a surjection. Give $M$ a basepoint that is disjoint from $N(S_1) \cup N(S_2)$.
We just have to find based loops $\ell_1$ and $\ell_2$ in $M$ such that
the image of $\ell_i$ is the generator for the $i$th factor of $\pi_1(({\Bbb R}P^3) \vee ({\Bbb R}P^3))$.
Pick a path from the basepoint of $M$ to one of the 0-cells of $N(S_i)$, so that the interior of the
path misses $N(S_1) \cup N(S_2)$. Then continue this path across the 1-cell that intersects $S_i$. Then run the
path back to the basepoint of $M$, again with the interior of the path avoiding $N(S_1) \cup N(S_2)$.
This is possible because we are assuming that $M - (S_1 \cup S_2)$ is connected.
The result is a based loop $\ell_i$ with the required properties.

(i) $\Rightarrow$ (ii). Suppose that there is a surjective homomorphism 
$\phi \colon \pi_1(M) \rightarrow {\Bbb Z}_2 \ast {\Bbb Z}_2$.  

Consider ${\Bbb R}P^2 \vee {\Bbb R}P^2$.
Give each copy of ${\Bbb R}P^2$ the standard
cell structure, and suppose that these two copies of ${\Bbb R}P^2$
are glued along their 0-cells. In the interior of each 1-cell, pick a point, and in each
2-cell, pick a properly embedded arc with endpoints equal to
one of these points. Let $\alpha_1$ and $\alpha_2$ be the
resulting disjoint simple closed curves in ${\Bbb R}P^2 \vee {\Bbb R}P^2$.

Pick a triangulation $T$ for $M$, and let $T^{(1)}$ and $T^{(2)}$ be its
1-skeleton and 2-skeleton. It is shown in the proof of Theorem 3.6 in [3]
that $\phi$ is induced by a map $f \colon T^{(2)} \rightarrow
{\Bbb R}P^2 \vee {\Bbb R}P^2$ with the following properties:
\item{(1)} For $i = 1$ and $2$, $f^{-1}(\alpha_i)$ is disjoint from the 0-skeleton
of $T$ and intersects the 1-skeleton in finitely many points.
Moreover, $f^{-1}(\alpha_i)$ intersects each face
of $T$ in a collection of properly embedded arcs,
with boundary equal to $f^{-1}(\alpha_i) \cap T^{(1)}$.
(In [3], this arrangement is called a regular mod 2 cocycle,
where each interior vertex has valence 2. Note also the space
referred to as $L(2)$ in [3] is just ${\Bbb R}P^2$.)
\item{(2)} $f^{-1}(\alpha_1 \cup \alpha_2)$ is non-separating in $T^{(2)}$.

We now extend $f^{-1}(\alpha_1)$ and $f^{-1}(\alpha_2)$
to disjoint surfaces $S_1$ and $S_2$ properly embedded in $M$.
For each tetrahedron $\Delta$ of $T$, $f^{-1}(\alpha_1)$ and
$f^{-1}(\alpha_2)$ intersect this tetrahedron in a collection
of simple closed curves in the boundary of $\Delta$.
We attach a collection of disjoint discs properly embedded
in $\Delta$ to these curves.

Note that $S_1 \cup S_2$ is non-separating. For consider
two points in the complement of $S_1 \cup S_2$.
We may find paths in $M - (S_1 \cup S_2)$
from these points to the 2-skeleton of $T$.
Because $f^{-1}(\alpha_1 \cup \alpha_2)$
is non-separating in $T^{(2)}$, we may find a path
joining these two points in the complement of $S_1 \cup S_2$.
$\square$

\noindent {\bf Remark 3.2.} Note that we can gain control over the mod 2 homology
classes of the surfaces $S_1$ and $S_2$ in Theorem 3.1, in terms of the
surjective homomorphism $\phi \colon \pi_1(M) \rightarrow {\Bbb Z}_2 \ast {\Bbb Z}_2$.
We see from the construction that the composition of $\phi$ with
projection ${\Bbb Z}_2 \ast {\Bbb Z}_2 \rightarrow {\Bbb Z}_2$
onto the $i$th factor is equal to the homomorphism $\pi_1(M) \rightarrow
{\Bbb Z}_2$ that counts the mod 2 intersection number with $S_i$.

We can now prove one direction of Theorem 1.1. Suppose that $M$ is a compact connected
orientable 3-manifold and that $\partial M$ is non-empty and contains no 2-spheres.
Suppose also that $M$ is neither a ${\Bbb Z}_2$ homology solid torus nor a ${\Bbb Z}_2$
homology cobordism between two tori. Then, we wish to show that $M$ contains 
two disjoint properly embedded surfaces $S_1$ and $S_2$ in $M$
such that $M - (S_1 \cup S_2)$ is connected. By Theorem 3.1, this 
is equivalent to $\pi_1(M)$ admitting a surjective homomorphism onto
${\Bbb Z}_2 \ast {\Bbb Z}_2$. By Corollary 2.2, this is equivalent
to the existence of an index 2 subgroup $K$ of $\pi_1(M)$ such that
$b_1(K) > b_1(M)$. To find such a subgroup, we use the following lemma.

\noindent {\bf Lemma 3.3.} {\sl Let $M$ be a compact orientable 3-manifold,
and let $\tilde M$ be a double cover of $M$ such that ${\rm genus}(\partial \tilde M) > {\rm genus}(\partial M)$.
Then $b_1(\tilde M) > b_1(M)$.}

\noindent {\sl Proof.} Let $p \colon \tilde M \rightarrow M$ be the covering map. It is argued
in the proof of Theorem 2.1 that $p^\ast \colon H^1(M; {\Bbb R})
\rightarrow H^1(\tilde M; {\Bbb R})$ is an injection. The same is
true of $(p|\partial \tilde M)^\ast \colon H^1(\partial M; {\Bbb R})
\rightarrow H^1(\partial \tilde M ; {\Bbb R})$. However, an alternative
proof is required, because $\partial M$ may be disconnected
and so an argument involving the fundamental group is not immediately
appropriate. Instead, we consider the transfer homomorphism
$t \colon H^1(\partial \tilde M; {\Bbb R}) \rightarrow H^1(\partial M; {\Bbb R})$.
Recall that this defined by sending a cocycle $c \in C^1(\partial \tilde M; {\Bbb R})$
to the cocycle $tc \in C^1(\partial M; {\Bbb R})$, where $tc(e) = c(\tilde e_1) + c(\tilde e_2)$.
Here $e$ is an oriented edge in some cell structure on $\partial M$,
and $e_1$ and $e_2$ are its inverse images in $\partial \tilde M$.
It is well known that this gives a well-defined homomorphism
$t \colon H^1(\partial \tilde M; {\Bbb R}) \rightarrow H^1(\partial M; {\Bbb R})$.
It is clear the composition $t \circ (p|\partial \tilde M)^\ast \colon H^1(\partial M; {\Bbb R})
\rightarrow H^1(\partial M; {\Bbb R})$ is the homomorphism that multiplies
by $2$. This is implies that $(p|\partial \tilde M)^\ast$ is an injection, as required.

Now consider the commutative diagram
$$\matrix{
H^1(M; {\Bbb R}) &\buildrel p^\ast \over \longrightarrow &H^1(\tilde M; {\Bbb R}) \cr
\downarrow && \downarrow \cr
H^1(\partial M; {\Bbb R}) &\buildrel (p|\partial \tilde M)^\ast \over \longrightarrow 
&H^1(\partial \tilde M; {\Bbb R})}
$$
where the vertical arrows are the homomorphisms induced by inclusion.
It is a well-known consequence of Poincar\'e duality
that, for the compact orientable 3-manifold $M$, the image of
$H^1(M; {\Bbb R}) \rightarrow H^1(\partial M; {\Bbb R})$ has dimension equal to exactly
half the dimension of $H^1(\partial M; {\Bbb R})$. A similar statement is
true for $\tilde M$.

Suppose now that $b_1(\tilde M) \leq b_1(M)$. Then
$p^\ast \colon H^1(M; {\Bbb R}) \rightarrow H^1(\tilde M; {\Bbb R})$
is therefore an isomorphism. So,
$$\eqalign{
&{\rm Im}(H^1(\tilde M; {\Bbb R}) \rightarrow H^1(\partial \tilde M;{\Bbb R})) \cr
= \  &{\rm Im}(H^1(M; {\Bbb R}) \rightarrow H^1(\partial \tilde M;{\Bbb R})) \cr
\cong \ & {\rm Im}(H^1(M; {\Bbb R}) \rightarrow H^1(\partial M;{\Bbb R})). \cr}$$
Hence, $b_1(\partial \tilde M) = b_1(\partial M)$. But this
is equivalent to the statement that $\partial \tilde M$ and $\partial M$
have the same genus, which is contrary to hypothesis. $\square$

We now return to the proof of one direction of Theorem 1.1.

\noindent {\sl Case 1.} $\partial M$ has a component which is not a torus.

Then, $\chi(\partial M) < 0$. Pick any double cover $\tilde M \rightarrow M$.
Then, $\chi(\partial \tilde M) = 2 \chi(\partial M) < \chi(\partial M)$.
Also, $|\partial \tilde M| \geq |\partial M|$. This implies that 
$${\rm genus}(\partial \tilde M)  = {- \chi(\partial \tilde M) + 2 |\partial \tilde M| \over 2}
> {- \chi(\partial M) + 2 |\partial M| \over 2} = {\rm genus}(\partial M),$$
as required.

\noindent {\sl Case 2.} $\partial M$ is at least three tori.

Then, $H^1(M; {\Bbb Z}_2)$ has rank at least three. Pick some torus $T$ in $\partial M$.
Then the homomorphism $H^1(M; {\Bbb Z}_2) \rightarrow H^1(T; {\Bbb Z}_2)$
that is induced by inclusion has non-zero kernel. Pick some non-zero
element in this kernel, and let $\tilde M$ be the corresponding double
cover of $M$. Then, the inverse image of $T$ in $\tilde M$ is two copies of
$T$. So, ${\rm genus}(\partial \tilde M) = |\partial \tilde M| > |\partial M| > {\rm genus}(\partial M)$, as
required.

\noindent {\sl Case 3.} $\partial M$ consists of two tori.

So, $H^1(M; {\Bbb Z}_2)$ has rank at least two. Now it cannot be
the case that, for each torus $T$ in $\partial M$, $H^1(M; {\Bbb Z}_2)
\rightarrow H^1(T; {\Bbb Z}_2)$ is an isomorphism. For this would
imply that $M$ is a ${\Bbb Z}_2$ homology cobordism between the two components
of $\partial M$, and this is contrary to assumption. So, for some
component $T$ of $\partial M$,
$H^1(M; {\Bbb Z}_2) \rightarrow H^1(T; {\Bbb Z}_2)$ is not an isomorphism,
and is therefore not injective. As in Case 2, we consider a double
cover $\tilde M$ of $M$ corresponding to a non-zero element in the kernel.
This has ${\rm genus}(\partial \tilde M) > {\rm genus}(\partial M)$, as
required.

\noindent {\sl Case 4.} $\partial M$ is a single torus $T$.

Then, $H^1(M; {\Bbb Z}_2)$ has rank at least $1$. In fact, it
must have rank at least $2$, since otherwise $M$ is a ${\Bbb Z}_2$ homology
solid torus. But the image of $H^1(M; {\Bbb Z}_2) \rightarrow H^1(T; {\Bbb Z}_2)$
has rank one, and so again, there is a non-trivial element in its
kernel. The argument then proceeds as in Cases 2 and 3. $\square$

\vskip 18pt
\centerline{\caps 4. Making the construction explicit}
\vskip 6pt

In the previous section, we completed the proof of one direction of Theorem 1.1, thereby
establishing the existence of pairs of disjoint surfaces in many 3-manifolds whose union is non-separating.
However, the proof is rather algebraic, and so it is hard to see 
how the surfaces arise explicitly. In this section, we will remedy this
defect, by providing an alternative way of constructing these surfaces
that is considerably more geometric.

A key part of the construction was to find a double cover $\tilde M \rightarrow M$
for which $b_1(\tilde M) > b_1(M)$. We saw in the proof of Theorem 2.1
that there is then a non-trivial primitive element $\alpha$ of $H^1(\tilde M; {\Bbb Z})$ that is in the
$-1$ eigenspace of the action of the non-trivial covering transformation $\tau$.

We will show that the surfaces $S_1$ and $S_2$ required by Theorem 1.1
may be constructed as follows:
\item{1.} Find a compact, oriented, properly embedded, non-separating
surface $S$ in $\tilde M$ that is dual to $\alpha$ and that is
invariant under the covering transformation $\tau$, but for which $\tau$ reverses the orientation.
We will prove below that such a surface $S$ always exists. Since $\tau$ reverses the
orientation of $S$, its image in $M$ is unoriented and typically non-orientable. This image will be
one of the surfaces $S_1$.
\item{2.} Cut $\tilde M$ along $S$ to give a compact orientable 3-manifold $M'$. The involution
$\tau$ restricts to $M'$, and it swaps the two copies of $S$ in $\partial M'$ (called $S_-$ and $S_+$, say). 
The next stage is to find a properly embedded surface $S'$ in $M'$ that is disjoint from $S_- \cup S_+$ 
and that is invariant under $\tau$. It must separate $M'$ into two components,
one containing $S_-$, the other containing $S_+$. Thus, $\tau$ swaps these
two components. Again, the existence of such a surface
will be established below. In practice, it is not hard to find. The image of $S'$ in $M$
will be $S_2$.

Note that if these surfaces $S$ and $S'$ exist, as claimed above, then their images
$S_1$ and $S_2$ in $M$ have the required properties. Note that $S$ and $S'$
are properly embedded and disjoint, by construction, and are the inverse
images of $S_1$ and $S_2$. Hence, $S_1$ and $S_2$ are also properly embedded
and disjoint. Also, $\tilde M - (S \cup S')$ has two components that
are swapped by $\tau$. So, $M - (S_1 \cup S_2)$ is connected, as required.

We now show that the surfaces $S$ and $S'$ always exist. We know from Theorem 1.1
that $M$ contains properly embedded disjoint surfaces $S_1$ and $S_2$ such
that $M - (S_1 \cup S_2)$ is connected. Also, from the proof of Theorem 1.1,
$\tilde M$ is the double cover of $M$ corresponding to $[S_1] + [S_2] \in H_2(M, \partial M; {\Bbb Z}_2)$.
This cover is constructed as follows. Let $\tilde S_1$ and $\tilde S_2$ be
${\rm cl}(\partial N(S_1) - \partial M)$ and ${\rm cl}(\partial N(S_2) - \partial M)$.
These are (possibly disconnected) double covers of $S_1$ and $S_2$,
respectively. Then $\tilde M$ is constructed by gluing together two copies of 
${\rm cl}(M - N(S_1 \cup S_2))$,
a copy of $\tilde S_1 \times [-1,1]$ and a copy of $\tilde S_2 \times [-1,1]$, as follows.
We attach $\tilde S_1 \times \{1\}$ and $\tilde S_2 \times \{1\}$ to one copy of
${\rm cl}(M - N(S_1 \cup S_2))$, using the identity map. But we attach 
$\tilde S_1 \times \{ -1 \}$ to the
other copy of ${\rm cl}(M - N(S_1 \cup S_2))$ via the covering
involution on $\tilde S_1$. We attach $\tilde S_2 \times \{ -1\}$
in a similar way.

We take $S$ and $S'$ to be the inverse image in $\tilde M$ of $S_1$ and $S_2$.
These have the required properties. For example, $S$ is $\tilde S_1 \times \{ 0 \}$
which is transversely oriented in the product bundle $\tilde S_1 \times [-1,1]$.
Hence it is oriented. Also, the covering transformation on $\tilde M$ reverses
this transverse orientation, and hence reverses the orientation on $S$.
Note that $S$ is an oriented surface, properly embedded and non-separating
in $\tilde M$. Hence, it represents a non-trivial primitive element of $H^1(\tilde M; {\Bbb Z})$.
Since $\tau$ preserves $S$ but reverses its orientation, this class is in
the $-1$ eigenspace.

Thus, the existence of $S$ and $S'$ is proved using the existence of $S_1$ and $S_2$.
But in practice, it is easiest to find $S$ and $S'$ first, and from these, construct $S_1$ and $S_2$.

We give a concrete example. Let $M$ be the exterior of the Whitehead link $L_1 \cup L_2$, shown in the left in Figure 1. 
Let $\tilde M$ be the double cover of $M$ that corresponds to the kernel of the
homomorphism $\pi_1(M) \rightarrow {\Bbb Z}_2$ which counts linking number mod $2$
with $L_2$. Since $L_2$ is unknotted, $\tilde M$ is the exterior of link in $S^3$
shown in the right of Figure 1. The link has three components, whereas the Whitehead link
has two, and so $b_1(\tilde M) > b_1(M)$. The surface $S$ shown in the right of Figure 1 is
non-separating, orientable and properly embedded in $\tilde M$. The covering involution
$\tau$ of $\tilde M$ preserves $S$ but reverses its orientation. Thus, Step 1 above
applies, and we may take $S_1$ to be the image of $S$ in $M$. This is shown in the left of Figure 1.

\vskip 18pt
\centerline{
\includegraphics[width=3.5in]{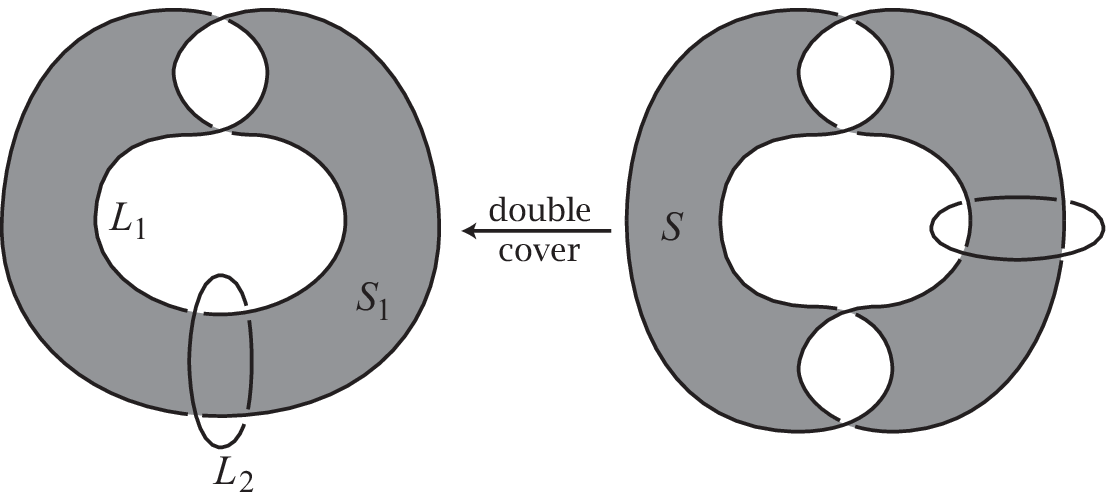}
}
\centerline{Figure 1.}

The surface $S'$ is not quite so easy to see. Let $D$ be the disc properly embedded in
the exterior of $S_1$ shown in Figure 2. Its inverse image in $\tilde M$ is two discs
$\tilde D$ properly embedded in the exterior of $S$. Now $M'$ (the exterior of $S$)
is a sutured manifold because when $S$ is oriented,
the two copies of $S$ in $\partial M'$ naturally point into and out of
$M'$. Let $\gamma'$ be its sutures. The discs $\tilde D$ form product discs. 
We orient them in such a way that this orientation is reversed by the
covering involution. Let $M_2'$ be obtained
by decomposing $M'$ along $\tilde D$. This is homeomorphic to $T^2 \times [0,1]$.
One component of its boundary, say $T^2 \times \{0 \}$, contains two sutures $\gamma_2'$;
the other has none. The covering involution $\tau$ restricts to an involution of
$M_2'$, which preserves the product structure on $T^2 \times [0,1]$ and preserves the sutures.
However, it swaps the inward and outward-pointing parts of $\partial M_2'$.
Let $A'$ be the vertical annuli $\gamma_2' \times [0,1]$ in $T^2 \times [0,1]$. 
Now, inside $\partial M_2'$ lies two copies of $\tilde D$, which therefore forms
four discs. The sutures $\gamma_2'$ run over each of these discs in a single arc.

\vskip 18pt
\centerline{
\includegraphics[width=3.5in]{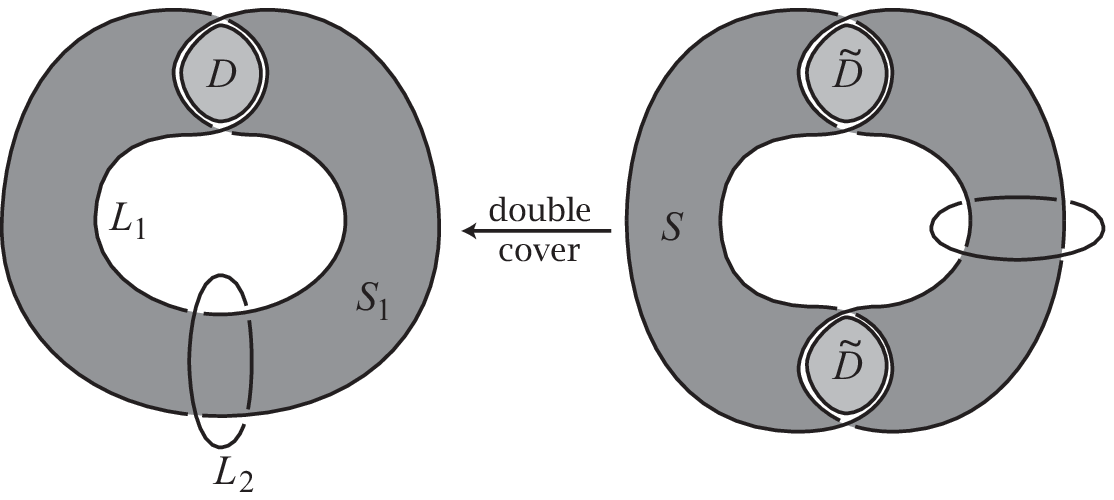}
}
\centerline{Figure 2.}

When we reverse the sutured manifold decomposition, to form
$M'$ from $M_2'$, we attach $\tilde D \times I$ to $M_2'$.
Inside each component of $\tilde D \times I$, we may insert a band of the
form $[0,1] \times[0,1]$ where $[0,1] \times \{0, 1\}$ lies in $\gamma_2'$
and $\{0,1\} \times [0,1]$ lies in $\gamma'$ (see Figure 3). Attaching these two
strips to $A'$ forms the required surface $S'$. It is properly embedded
in $M'$ and is disjoint from $S$. It is invariant under the action of
$\tau$. Also, it separates $M'$ into two components, which are swapped by $\tau$.

\vskip 18pt
\centerline{
\includegraphics[width=3.5in]{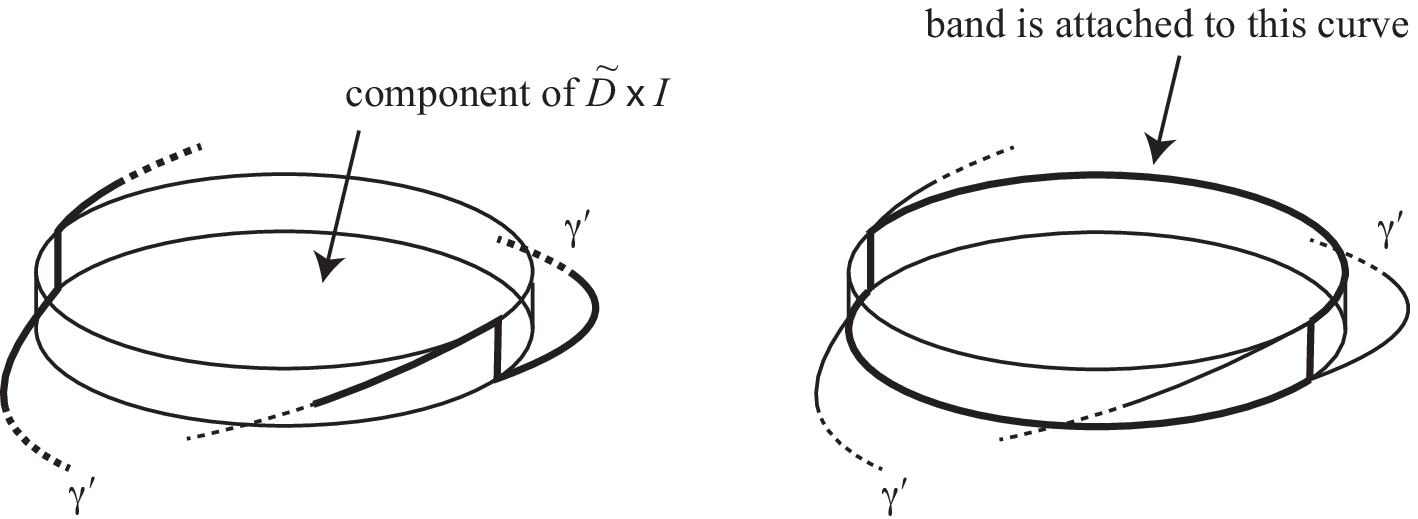}
}
\centerline{Figure 3.}

Let $S_2$ be the image of $S'$ in $M$. This can be seen as follows. 
The image of $\tilde D$ in the exterior of $S_1$ is the disc $D$. One cannot call it a product
disc because the exterior of $S_1$ is not a sutured manifold. Nevertheless, one can cut
the exterior of $S_1$ along $D$ to form a space $M_2$ which is also homeomorphic to $T^2 \times I$.
This is the quotient of $M_2'$ under $\tau$. The image of $\gamma_2'$
is a single curve in $\partial M_2$. Let $A$ be the vertical annulus over this in $T^2 \times I$. This
is properly embedded in $M_2$. Reconstruct the exterior of $S_1$ from $M_2$ by
reattaching $D \times I$. Inside $D \times I$, we may find a band.
Attaching this band to $A$ gives the surface $S_2$. It is a twice-punctured
projective plane. It has two boundary components, one lying in $L_1$
and one lying in $L_2$. It is disjoint from $S_1$, and $S_1 \cup S_2$
is non-separating.

The surface $S_1$ is easy to see, but the difficulty in visualising $S_2$
perhaps arises from the fact that its boundary component on $L_2$
has slope $2/1$. Note that it is a spanning surface for $L_1 \cup L_2$.
In fact, it is shown in Figure 4, after an isotopy. This was constructed by
retracting the annulus $A$ a little so that it lies in a small regular
neighbourhood of $L_2$, and then attaching the band to form $S_2$.

\vskip 18pt
\centerline{
\includegraphics[width=1.5in]{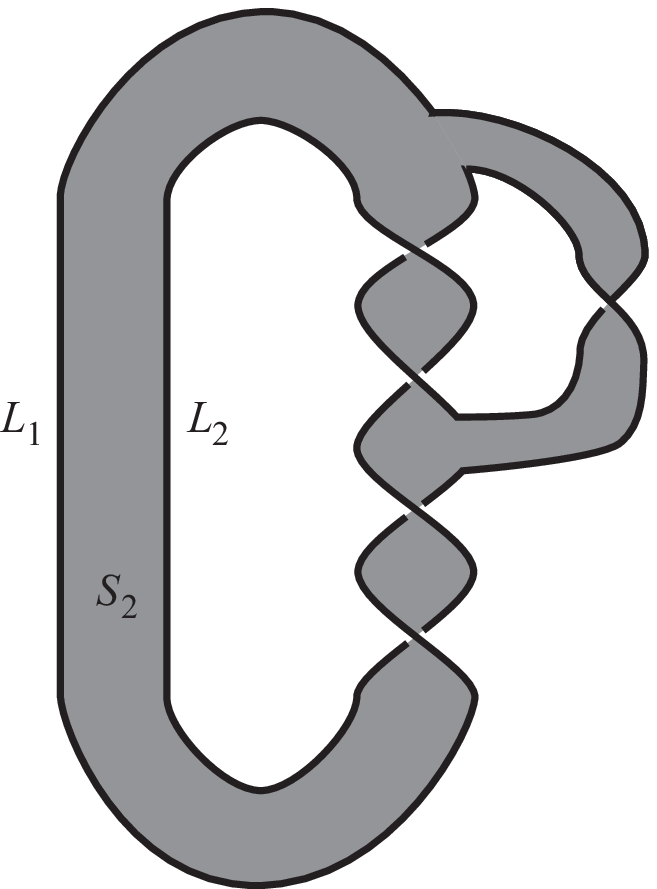}
}
\vskip 6pt
\centerline{Figure 4.}

\vskip 18pt 
\centerline{\caps 5. Profinite methods}
\vskip 6pt

We still need to prove one direction of Theorem 1.1. We will show that if $M$ is a compact orientable 3-manifold 
that is either a ${\Bbb Z}_2$ homology solid torus or a ${\Bbb Z}_2$ homology cobordism between two tori, 
then it cannot support the surfaces
$S_1$ and $S_2$ that are described in the theorem. In the course of this proof, we will introduce some 
techniques from the theory of pro-$p$ groups. These will turn out to have other uses. In particular, we will
be able to use them to gain control of the ${\Bbb Z}_2$ homology classes of the surfaces $S_1$ and $S_2$
for certain 3-manifolds $M$.

Let $p$ be a prime. Recall that the {\sl pro-$p$ completion} $\widehat \Gamma_{(p)}$ of a group $\Gamma$ is the inverse limit of all its
finite quotients that are $p$-groups. More precisely, an element of $\widehat \Gamma_{(p)}$ is a choice, for
each normal subgroup $N$ of $\Gamma$ with index a power of $p$, of an element $g_N$ of $\Gamma/N$,
subject to the following compatibility condition. Whenever $N$ and $N'$ are normal subgroups of $\Gamma$
with index that are powers of $p$, and satisfying $N \geq N'$, then we insist that $g_{N'}$ maps to $g_N$
under the quotient map $\Gamma/N' \rightarrow \Gamma/N$.

The definition of $\widehat \Gamma_{(p)}$ is phrased in terms of normal subgroups of $\Gamma$ with index
a power of $p$. However, $\widehat \Gamma_{(p)}$ contains information about a wider class of subgroups of $\Gamma$
which are defined as follows. We say that a subgroup $N$ of $\Gamma$ is {\sl co-$p$} if there is 
a sequence of finite index subgroups
$$\Gamma = N_0 \geq N_1 \geq \dots \geq N_k = N$$
such that each $N_i$ is normal in $N_{i-1}$ and has index a power of $p$. The terminology co-$p$ is not
standard. The usual phrase is `subnormal with index a power of $p$'. 

The following is a rapid consequence of the definition of $\widehat \Gamma_{(p)}$.

\noindent {\bf Proposition 5.1.} {\sl Let $G$ and $\Gamma$ be finitely generated discrete groups, and
let $p$ be a prime. Suppose that there is a group isomorphism $\phi \colon \widehat G_{(p)} \rightarrow \widehat \Gamma_{(p)}$
between their pro-$p$ completions. Then the following hold.

\item{(i)} There is an induced bijection (also denoted $\phi$) between the set of co-$p$ subgroups of $G$
and the set of co-$p$ subgroups of $\Gamma$.
\item{(ii)} For a co-$p$ subgroup $N$ of $G$, $[G:N] = [\Gamma: \phi(N)]$. 
\item{(iii)} If $N$ is any co-$p$ subgroup of $G$, then $N$ is normal in $G$ if and only if $\phi(N)$
is normal in $\Gamma$. In this case, $G/N$ is isomorphic to $\Gamma / \phi(N)$.
\item{(iv)} If $N$ and $N'$ are co-$p$ subgroups of $G$, then $N \subset N'$ if and only if $\phi(N) \subset \phi(N')$.
\item{(v)} If $N$ is a co-$p$ subgroup of $G$, then $N$ and $\phi(N)$ have isomorphic pro-$p$ completions.

}

Thus, an isomorphism between $\widehat G_{(p)}$ and $\widehat \Gamma_{(p)}$ entails a strong correspondence between the
co-$p$ subgroups of $G$ and the co-$p$ subgroups of $\Gamma$. Slightly surprisingly, $\widehat \Gamma_{(p)}$
also controls the first Betti number of $\Gamma$.

\noindent {\bf Proposition 5.2.} {\sl Let $G$ and $\Gamma$ be finitely generated groups with isomorphic
pro-$p$ completions for some prime $p$. Then $b_1(G) = b_1(\Gamma)$.}

\noindent {\sl Proof.} By (iii) of Proposition 5.1, there is a one-one correspondence between
the quotients of $G$ that are abelian $p$-groups
and the similar set of quotients of $\Gamma$. But $G$ surjects onto $({\Bbb Z}/p^k {\Bbb Z})^l$
for all $k \in {\Bbb N}$ if and only if $l \leq b_1(G)$. So, $b_1(G) = b_1(\Gamma)$.
$\square$

These results imply that, for a finitely generated group $\Gamma$, the existence of a surjective
homomorphism from $\Gamma$ onto the infinite dihedral group is determined by $\widehat \Gamma_{(2)}$.

\noindent {\bf Theorem 5.3.} {\sl Let $G$ and $\Gamma$ be finitely presented groups with isomorphic
pro-$2$ completions. Then $G$ admits a surjective homomorphism onto ${\Bbb Z}_2 \ast {\Bbb Z}_2$
if and only if $\Gamma$ does.}

\noindent {\sl Proof.} Suppose that $G$ admits a surjective homomorphism onto ${\Bbb Z}_2 \ast {\Bbb Z}_2$.
Then, by Corollary 2.2, $G$ has an index 2 subgroup $G'$ with $b_1(G') > b_1(G)$. Index 2
subgroups are normal and hence co-2. So, using Proposition 5.1, $G'$ corresponds to a subgroup $\Gamma'$ of
$\Gamma$ with index 2. Also by Proposition 5.1, $G'$ and $\Gamma'$ have isomorphic pro-$2$
completions. So, by Proposition 5.2, $b_1(G') = b_1(\Gamma')$. Similarly, $b_1(G) = b_1(\Gamma)$,
and hence $b_1(\Gamma') > b_1(\Gamma)$.
So, by Corollary 2.2, $\Gamma$ admits a surjective homomorphism onto ${\Bbb Z}_2 \ast {\Bbb Z}_2$.
$\square$

The following result shows that the existence of an isomorphism between pro-$p$ completions is surprisingly common.

\noindent {\bf Theorem 5.4.} {\sl Suppose that there is a homomorphism $\phi \colon G \rightarrow \Gamma$ between finitely
presented groups that induces an isomorphism $H_1(G; {\Bbb Z}_p) \rightarrow H_1(\Gamma; {\Bbb Z}_p)$
and a surjection $H_2(G ; {\Bbb Z}_p) \rightarrow H_2(\Gamma; {\Bbb Z}_p)$. Then 
$\phi$ induces an isomorphism between the pro-$p$ completions of
$G$ and $\Gamma$.}

This is a fairly well known consequence of work of Stallings [5], and a proof can be found in [4]
(see Theorem 2.12 of [4]).

Note that for any finite cell complex $M$, there is a surjective homomorphism 
$H_2(M; {\Bbb Z}_p) \rightarrow H_2(\pi_1(M); {\Bbb Z}_p)$. This is because 
$H_2(\pi_1(M); {\Bbb Z}_p)$ is the homology of an Eilenberg-Maclane space
$K(\pi_1(M); 1)$, which can be obtained from $M$ by attaching cells in dimensions
3 and higher.

We are now in a position to prove the remaining direction of Theorem 1.1. Let $M$ be a compact orientable
3-manifold with non-empty boundary that contains no 2-spheres.

Suppose first that $M$ is  a ${\Bbb Z}_2$ homology solid torus. In other words,
$H_1(M; {\Bbb Z}_2) = {\Bbb Z}_2$ and $H_2(M; {\Bbb Z}_2) = 0$. Let $\Gamma$ be $\pi_1(M)$.
Then $H_1(\Gamma; {\Bbb Z}_2) = {\Bbb Z}_2$ and $H_2(\Gamma; {\Bbb Z}_2) = 0$. Let
$\phi \colon {\Bbb Z} \rightarrow \Gamma$ be a homomorphism that sends a generator of ${\Bbb Z}$
to an element of $\Gamma$ that is non-trivial in $H_1(\Gamma; {\Bbb Z}_2)$. Then $\phi$ induces
isomorphisms on first and second homology with ${\Bbb Z}_2$ coefficients. Therefore,
by Theorem 5.4 the pro-$2$ completions of ${\Bbb Z}$ and $\Gamma$ are isomorphic.
Since there is no surjective homomorphism from ${\Bbb Z}$ to ${\Bbb Z}_2 \ast {\Bbb Z}_2$,
the same is therefore true for $\Gamma$, by Theorem 5.3.

Now suppose that $M$ is a ${\Bbb Z}_2$ homology cobordism between two tori $T_1$
and $T_2$. Then, by assumption, $i_\ast \colon H_1(T_1; {\Bbb Z}_2) \rightarrow H_1(M; {\Bbb Z}_2)$
is an isomorphism, where $i \colon T_1 \rightarrow M$ is inclusion.
We will also show that $i_\ast \colon H_2(T_1; {\Bbb Z}_2) \rightarrow H_2(M; {\Bbb Z}_2)$
is a surjection. Now there is an exact sequence
$$H_2(T_1; {\Bbb Z}_2) \buildrel i_\ast \over \rightarrow H_2(M; {\Bbb Z}_2) \rightarrow H_2(M,T_1; {\Bbb Z}_2)$$
and so it suffices to show that $H_2(M,T_1; {\Bbb Z}_2)$ is trivial. But this is isomorphic to
$H^1(M, T_2; {\Bbb Z}_2)$ by Poincar\'e duality. This fits into an exact sequence
$$H^0(M; {\Bbb Z}_2) \buildrel \cong \over \rightarrow H^0(T_2; {\Bbb Z}_2) \rightarrow
H^1(M, T_2; {\Bbb Z}_2) \rightarrow H^1(M; {\Bbb Z}_2) \buildrel \cong \over \rightarrow
H^1(T_2; {\Bbb Z}_2).$$
So, $H^1(M, T_2; {\Bbb Z}_2)$ is trivial, and hence we have shown that
$i_\ast \colon H_2(T_1; {\Bbb Z}_2) \rightarrow H_2(M; {\Bbb Z}_2)$ is a surjection.
Now the torus is an Eilenberg-Maclane space and so $H_2(T_1; {\Bbb Z}_2) = H_2(\pi_1(T_1); {\Bbb Z}_2)$.
Therefore, $i_\ast \colon H_2(\pi_1(T_1); {\Bbb Z}_2) \rightarrow H_2(\pi_1(M); {\Bbb Z}_2)$ is also a surjection.

Hence, using Theorem 5.4, $\pi_1(T_1)$ and $\pi_1(M)$ have isomorphic pro-$2$ completions.
Now, $\pi_1(T_1) = {\Bbb Z} \times {\Bbb Z}$ clearly does not admit a surjective homomorphism onto
${\Bbb Z}_2 \ast {\Bbb Z}_2$. Hence, by Theorem 5.3, nor does $\pi_1(M)$.

This completes the proof of Theorem 1.1. $\square$

\vskip 18pt
\centerline{\caps 6. Controlling the homology classes of the surfaces}
\vskip 6pt

Given how common it is for a compact orientable 3-manifold $M$ to contain disjoint properly embedded
surfaces $S_1$ and $S_2$ such that $M - (S_1 \cup S_2)$ is connected, it is natural
to ask which pairs of classes in $H_2(M, \partial M; {\Bbb Z}_2)$ may be represented by
such surfaces. In this section, we will address this question.

We start by observing that $[S_1]$ and $[S_2]$ must be non-trivial and distinct.
This because, when these surfaces exist, there is an associated surjective homomorphism
$\pi_1(M) \rightarrow {\Bbb Z}_2 \ast {\Bbb Z}_2$. Composing this with the surjection
onto the first and second factors, we obtain two homomorphisms $\phi_1$ and $\phi_2 \colon
\pi_1(M) \rightarrow {\Bbb Z}_2$. These correspond to two classes $[\phi_1]$ and
$[\phi_2]$ in $H^1(M; {\Bbb Z}_2)$ and the Poincar\'e duals of these classes
are $[S_1]$ and $[S_2]$ in $H_2(M, \partial M; {\Bbb Z}_2)$. Hence, because
$\phi_1$ and $\phi_2$ are distinct non-trivial homomorphisms to ${\Bbb Z}_2$,
we deduce that $[S_1]$ and $[S_2]$ are distinct and non-trivial.

How much control do we have over $[S_1]$ and $[S_2]$, given the above restriction?
Equivalently, how much control do we have over $\phi_1$ and $\phi_2$?
By Theorem 2.1, the kernel of $(\phi_1 + \phi_2) \colon \pi_1(M) \rightarrow {\Bbb Z}_2$ 
is an index two subgroup $K$ of
$\pi_1(M)$ such that $b_1(K) > b_1(M)$. Thus, we deduce the following
from Theorem 2.1, Theorem 3.1 and Remark 3.2.

\vfill\eject

\noindent {\bf Proposition 6.1.} {\sl Let $M$ be a compact connected orientable 3-manifold,
and let $\alpha$ be a non-trivial class in $H^1(M ; {\Bbb Z}_2)$.
Then, the following are equivalent:
\item{(i)} There are disjoint properly embedded surfaces $S_1$ and $S_2$
such that $M - (S_1 \cup S_2)$ is connected, and such that
$S_1 \cup S_2$ is dual to $\alpha$.
\item{(ii)} The degree two cover $\tilde M$ of $M$ that corresponds to
${\rm ker}(\alpha \colon \pi_1(M) \rightarrow {\Bbb Z}_2)$ satisfies
$b_1( \tilde M) > b_1(M)$.

}

This controls only the class $[S_1] + [S_2]$ in $H_2(M, \partial M; {\Bbb Z}_2)$.
But nevertheless, it does constrain which classes can arise, as in the
following theorem. The version of this for links of 2-components was
proved by Hillman (see Theorem 7.7 in [2]).

\noindent {\bf Theorem 6.2.} {\sl Let $M$ be the exterior of a link $L$ in the 3-sphere,
and let $\alpha \in H^1(M ; {\Bbb Z}_2)$ count the linking number mod 2 with $L$.
Then, the following are equivalent:
\item{(i)} There are disjoint properly embedded surfaces $S_1$ and $S_2$
such that $M - (S_1 \cup S_2)$ is connected, and such that
$S_1 \cup S_2$ is dual to $\alpha$.
\item{(ii)} $\Delta_L(-1,-1, \dots, -1) = 0$, where $\Delta_L(t_1, t_2, \dots, t_n)$ is the
Alexander polynomial of $L$.
\item{(iii)} $L$ has a disconnected compact spanning surface, no component of which
is closed.

}

\noindent {\sl Proof.} (i) $\Leftrightarrow$ (ii): It is well known that the Alexander polynomial of a link $L$
encodes information about the homology of the finite abelian covers of its exterior.
In particular, $\Delta_L(-1, \dots, -1)$ equals the determinant of a square presentation
matrix for $H_1(\widehat M)$, where $\widehat M$ is the double cover of
$S^3$ branched over $L$ (see p.121 of [2] for example). Hence, $b_1(\widehat M) > 0$ if and only if
$\Delta_L(-1, \dots, -1) = 0$. An elementary calculation gives that
$b_1(\widehat M) = b_1(\tilde M) - b_1(M)$, where $\tilde M$
is the double cover of $M$ dual to $\alpha$. Thus, applying Proposition 6.1
completes the proof. 

(iii) $\Rightarrow$ (i): If $S$ is a compact spanning surface for $L$, then its
restriction to $M$ is a properly embedded surface dual to $\alpha$. If $S$ is
disconnected, then so too is the surface in $M$. If $S$ is a spanning
surface with no closed components, the surface in $M$ is non-separating.

(i) $\Rightarrow$ (iii): Let $S_1$ and $S_2$ be disjoint, properly embedded surfaces
in $M$ such that $S_1 \cup S_2$ is dual to $\alpha$, and such that
$M - (S_1 \cup S_2)$ is connected. We will show how to
modify $S_1$ and $S_2$, so that the number of components of $M - (S_1 \cup S_2)$
does not increase, and so that afterwards, $S_1 \cup S_2$ intersects each
component of $\partial M$ in a single simple closed curve. These modifications will not
change the homology classes of $S_1$ and $S_2$ and so each meridian of $L$
will still have non-empty intersection with $S_1 \cup S_2$. Also, $S_1$ and $S_2$
will remain non-empty, and each component of $S_1 \cup S_2$ will have non-empty boundary. Thus, we will be
able to extend $S_1 \cup S_2$ into $N(L)$ to form a disconnected compact spanning surface for $L$
with no closed components.

We may first assume that $\partial S_1 \cup \partial S_2$ is a collection of
essential simple closed curves on $\partial M$. For if some component of
$\partial S_1 \cup \partial S_2$ is inessential, we may find one that bounds
a disc with interior disjoint from $S_1 \cup S_2$. Attach this disc to
$S_1 \cup S_2$ and push it a little into the interior of $M$ to make
the surfaces properly embedded. This does not change the number of components
of $M - (S_1 \cup S_2)$. It preserves the properties of these surfaces, but reduces
the number of boundary components.

Thus, $\partial S_1 \cup \partial S_2$ divides each component of $\partial M$ into annuli.
Consider a meridian $\mu$ for some component of $L$. This is a simple closed curve
on a toral component $T$ of $\partial M$. Arrange for $\mu$ to intersect
$\partial S_1 \cup \partial S_2$ as few times as possible.
The evaluation of $\mu$ under $\alpha$ is 
$1$, and hence it intersects $\partial S_1 \cup \partial S_2$
an odd number of times. It cannot therefore meet
$\partial S_1$, then $\partial S_2$, then $\partial S_1$,
and so on in an alternating fashion. Thus, there are two successive intersections
which both lie in $\partial S_1$, say. If these lie in distinct components of $\partial S_1$,
then these two curves cobound an annulus in $T$. We 
attach this annulus to $S_1 \cup S_2$. This preserves the properties
of $S_1 \cup S_2$ given in (i), but again reduces its number of boundary
components. The only other possibility is that $\partial S_1 \cup \partial S_2$ intersects $T$ in
a single simple closed curve. Modifying this curve by a mod $2$
homology, we may assume that it intersects $\mu$ just once.
Since this applies on every component of $\partial M$, we may extend
$S_1 \cup S_2$ to the required spanning surface for $L$. $\square$

In the case of the Whitehead link $L$, $\Delta_L(t_1, t_2) = (t_1-1)(t_2 - 1)$,
and hence $\Delta_L(-1,-1) \not=0$. Therefore, although there exist two
disjoint properly embedded surfaces $S_1$ and $S_2$ such that $M - (S_1 \cup S_2)$
is connected, there do not exist two such surfaces that together form a
spanning surface for $L$.

We have already seen that the existence of a surjective homomorphism from 
a finitely generated group $G$ to ${\Bbb Z}_2 \ast {\Bbb Z}_2$ is controlled
by the pro-2 completion of $G$. We may refine this further, as follows.

Let $G$ and $\Gamma$ be finitely generated groups. 
Suppose that their {pro-$2$} completions are isomorphic.
By Proposition 5.1, this sets up an isomorphism between 
$G/([G,G]G^2)$ and $\Gamma/([\Gamma,\Gamma]\Gamma^2)$.
This is because $G/([G,G]G^2)$ may be characterised as the
quotient of $G$ that is an elementary abelian 2-group of
maximal rank, and there is a similar characterisation of 
$\Gamma/([\Gamma,\Gamma]\Gamma^2)$.
This therefore gives a 1-1 correspondence between
homomorphisms $G \rightarrow {\Bbb Z}_2$ and homomorphisms $\Gamma \rightarrow {\Bbb Z}_2$.

\noindent {\bf Proposition 6.3.} {\sl Let $G$ and $\Gamma$ be finitely
generated groups. Suppose that their {pro-$2$} completions are isomorphic.
Let $\phi_1$ and $\phi_2$ be non-trivial homomorphisms $G \rightarrow {\Bbb Z}_2$
and let $\phi_1'$ and $\phi_2'$ be the corresponding
homomorphisms $\Gamma \rightarrow {\Bbb Z}_2$. Suppose that there is a surjective
homomorphism $\phi \colon G \rightarrow {\Bbb Z}_2 \ast {\Bbb Z}_2$ such that the
composition with projection onto the $i$th factor is $\phi_i$, for $i = 1$ and $2$. Then
there is a surjective homomorphism $\Gamma \rightarrow {\Bbb Z}_2 \ast {\Bbb Z}_2$ such that the
composition with projection onto the $i$th factor is $\phi_i'$, for $i = 1$ and $2$.}

\noindent {\sl Proof.} We define two subgroups of $G$. Let $G_2$ be the
kernel of $G \buildrel \phi \over \rightarrow {\Bbb Z}_2 \ast {\Bbb Z}_2 
\buildrel \pi \over \rightarrow {\Bbb Z}_2$,
where $\pi$ sends the non-trivial element of each factor onto the non-trivial element
of ${\Bbb Z}_2$. This is an index 2 normal subgroup of $G$. Let $G_4$
be the kernel of $G \rightarrow {\Bbb Z}_2 \ast {\Bbb Z}_2 \rightarrow {\Bbb Z}_2 \times {\Bbb Z}_2$,
where again the first map is $\phi$, and the second map is abelianisation.
This is an index $4$ normal subgroup of $G$. It is equal to the elements of
$G$ that have trivial images under both $\phi_1$ and $\phi_2$.

By Proposition 5.1, $G_2$ (respectively, $G_4)$ corresponds to an index 2 (respectively, 4)
normal subgroup $\Gamma_2$ (respectively, $\Gamma_4$) of $\Gamma$. 
By Proposition 5.1 and 5.2, $b_1(G_2) = b_1(\Gamma_2)$, and therefore
$H^1(G_2; {\Bbb Z})$ and $H^1(\Gamma_2; {\Bbb Z})$ are isomorphic. However, we wish
to set up a slightly more precise correspondence between these two
cohomology groups.

The group $G/G_2$ acts on $H^1(G_2; {\Bbb R})$ by conjugation.
The non-trivial element of $G/G_2$ specifies an involution of  
$H^1(G_2; {\Bbb R})$, and hence $H^1(G_2; {\Bbb R})$ decomposes
into a direct sum of $+1$ and $-1$ eigenspaces $H^1(G_2; {\Bbb R})^+$
and $H^1(G_2; {\Bbb R})^-$. Let $H^1(G_2; {\Bbb Z})^-$ be the intersection
between $H^1(G_2; {\Bbb R})^-$ and the lattice $H^1(G_2; {\Bbb Z})$.
Define $H^1(\Gamma_2; {\Bbb Z})^-$ similarly. Our goal is to find an
isomorphism $H^1(G_2; {\Bbb Z})^- \rightarrow H^1(\Gamma_2; {\Bbb Z})^-$
such the following diagram commutes:

$$\matrix{
H^1(G_2; {\Bbb Z})^- &\rightarrow &H^1(\Gamma_2; {\Bbb Z})^- \cr
\downarrow && \downarrow \cr
H^1(G_2; {\Bbb Z}_2) &\rightarrow &H^1(\Gamma_2; {\Bbb Z}_2)}
$$
Here, the top horizontal arrow is the required isomorphism.
The vertical arrows are the composition of inclusion into the full integral
cohomology group, followed by reduction mod 2.
The bottom horizontal arrow is the isomorphism that arises
from the fact that the pro-2 completions of $G_2$ and $\Gamma_2$ are isomorphic.
The point is that $\phi | G_2$ gives a primitive element of $H^1(G_2; {\Bbb Z})^-$,
which we want to correspond to a primitive element of $H^1(\Gamma_2; {\Bbb Z})^-$.
This then determines a surjective homomorphism $\phi' \colon \Gamma \rightarrow {\Bbb Z}_2 \ast {\Bbb Z}_2$.
It is the commutativity of the above diagram that will ensure that the
composition onto the first and second factors will be $\phi_1'$ and $\phi_2'$.

Let $k>1$ be a large enough integer so that $G_2$ does not surject onto
$({\Bbb Z}/2^k {\Bbb Z})^{b_1(G_2) + 1}$, and similarly for $\Gamma_2$.
Let $\tilde G_2$ be the subgroup of $G_2$ generated by elements of $G_2$
that have finite order in $G_2 / [G_2,G_2]$, together with $2^k$th powers
in $G_2$. This is a normal subgroup of $G_2$, 
such that $G_2 / \tilde G_2$ is isomorphic to $({\Bbb Z}/2^k {\Bbb Z})^{b_1(G_2)}$.
Then, using the isomorphism between pro-2 completions,
$\tilde G_2$ corresponds to $\tilde \Gamma_2$, which is 
a normal subgroup of $\Gamma_2$
such that $\Gamma_2 / \tilde \Gamma_2$ is 
isomorphic to $({\Bbb Z}/2^k {\Bbb Z})^{b_1(\Gamma_2)}$.
Let $(G_2 / \tilde G_2)^\ast$ be the set of homomorphisms
from $G_2$ to ${\Bbb Z}/2^k{\Bbb Z}$. This is also isomorphic
to $({\Bbb Z}/2^k {\Bbb Z})^{b_1(G_2)}$. Define
$(\Gamma_2 / \tilde \Gamma_2)^\ast$ similarly.

Consider the groups $G/\tilde G_2$ and $\Gamma/\tilde \Gamma_2$.
These are isomorphic. So, the conjugation action of $G/G_2$
on $(G_2 / \tilde G_2)^\ast$ is equivalent to the conjugation action of 
$\Gamma/\Gamma_2$ on $(\Gamma_2 / \tilde \Gamma_2)^\ast$.
Let $(G_2 / \tilde G_2)^\ast_-$ be those elements of $(G_2 / \tilde G_2)^\ast$
that are sent to their inverses by this action. This is a subgroup
of $(G_2 / \tilde G_2)^\ast$. Define $(\Gamma_2 / \tilde \Gamma_2)^\ast_-$
similarly. Then we have an isomorphism between $(G_2 / \tilde G_2)^\ast_-$
and $(\Gamma_2 / \tilde \Gamma_2)^\ast_-$. Now, $(G_2 / \tilde G_2)^\ast$
can be identified with the reduction mod $2^k$ of $H^1(G_2; {\Bbb Z})$.
However, $(G_2 / \tilde G_2)^\ast_-$ is not necessarily the same as
the reduction mod $2^k$ of $H^1(G_2; {\Bbb Z})^-$. This is because
every element of $(G_2 / \tilde G_2)^\ast$ with order 2 is sent to
its inverse under the conjugation action. Therefore, let $(G_2 / \tilde G_2)^\ast_0$
be the set of elements in $(G_2 / \tilde G_2)^\ast_-$ that are
a multiple of an element of $(G_2 / \tilde G_2)^\ast_-$ which does
not have order 2. This is the image of $H^1(G_2; {\Bbb Z})^-$
under the mod $2^k$ reduction map. Define $(\Gamma_2 / \tilde \Gamma_2)^\ast_0$
similarly. Then $(G_2 / \tilde G_2)^\ast_0$ and $(\Gamma_2 / \tilde \Gamma_2)^\ast_0$
are isomorphic. Pick a basis for $H^1(G_2;{\Bbb Z})^-$ that maps to 
a minimal generating set for $(G_2 / \tilde G_2)^\ast_0$.
This is sent to a minimal generating set for $(\Gamma_2 / \tilde \Gamma_2)^\ast_0$.
We may pick elements in $H^1(\Gamma_2;{\Bbb Z})^-$
in their inverse image which form a basis for 
$H^1(\Gamma_2;{\Bbb Z})^-$. Mapping the basis elements of
$H^1(G_2;{\Bbb Z})^-$ to these basis elements of 
$H^1(\Gamma_2;{\Bbb Z})^-$ gives the required isomorphism
between $H^1(G_2;{\Bbb Z})^-$ and $H^1(\Gamma_2;{\Bbb Z})^-$.

Using this isomorphism, $\phi |G_2 \colon G_2 \rightarrow {\Bbb Z}$ corresponds
to a primitive element of $H^1(\Gamma_2; {\Bbb Z})^-$ and hence, as in the proof of Theorem 1.1,
we obtain a surjective homomorphism $\phi' \colon \Gamma \rightarrow {\Bbb Z}_2 \ast {\Bbb Z}_2$.

We now show that the compositions of $\phi'$ with projections onto the
first and second factors are $\phi_1'$ and $\phi_2'$.

Note that $G_4$ is the kernel of $G_2 \rightarrow {\Bbb Z} \rightarrow {\Bbb Z}_2$.
Here, the first homomorphism is the restriction of $\phi$ to $G_2$ and the second
is reduction mod $2$. Consider $\Gamma_2 \rightarrow {\Bbb Z} \rightarrow {\Bbb Z}_2$,
where the first homomorphism is the restriction of $\phi'$ to $\Gamma_2$ and the second
is reduction mod $2$. Due to the commutativity of the above diagram,
the kernel of this is precisely $\Gamma_4$.
Now, $\phi_1$ and $\phi_2$ can be characterised as the only non-trivial
homomorphisms $G \rightarrow {\Bbb Z}_2$ that are non-trivial on
$G_2$ but trivial on $G_4$.  A similar statement holds for $\phi'_1$
and $\phi'_2$. So, $\phi_1$ and $\phi_2$ do indeed correspond to
$\phi'_1$ and $\phi'_2$.

There is a minor complication. It may be the case that if we compose
$\phi' \colon G \rightarrow {\Bbb Z}_2 \ast {\Bbb Z}_2$ with projection onto
the first (respectively, second) factor then we get $\phi'_2$ (respectively, $\phi'_1$).
If this is the case, then redefine $\phi'$ by composing with the
automorphism of ${\Bbb Z}_2 \ast {\Bbb Z}_2$ that swaps the factors. $\square$

We will now show that, in an important special case, there is a great deal
of flexibility over the classes $[S_1]$ and $[S_2]$.

\noindent {\bf Theorem 6.4.} {\sl Let $M$ be a compact orientable
3-manifold with the same ${\Bbb Z}_2$ homology as a handlebody. Then, given any two distinct non-zero classes
$z_1$ and $z_2$ in $H_2(M, \partial M; {\Bbb Z}_2)$, there
are disjoint properly embedded surfaces $S_1$ and $S_2$ in $M$,
such that $M - (S_1 \cup S_2)$ is connected, and such that
$z_i = [S_i]$.}

\noindent {\sl Proof.} Pick a basepoint for $M$, and pick
based loops $\ell_1, \dots, \ell_n$ such that $[\ell_1], \dots, [\ell_n]$
form a basis for $H_1(M; {\Bbb Z}_2)$.
We may choose these loops so that the mod 2 intersection numbers
satisfy $[\ell_1].z_i = \delta_{1i}$, $[\ell_2].z_i = \delta_{2i}$
and $[\ell_j].z_i = 0$ for $j \geq 3$. Let $\phi_1'$ and $\phi_2' \colon \pi_1(M) \rightarrow {\Bbb Z}_2$
be the mod 2 intersection numbers with $z_1$ and $z_2$ respectively.

Let $F$ be the free group on $n$ generators. Let $\psi \colon F \rightarrow \pi_1(M)$
be the homomorphism that sends the $i$th free generator to $\ell_i$.
Then, $\psi$ induces an isomorphism $H_1(F ; {\Bbb Z}_2) \rightarrow H_1(\pi_1(M); {\Bbb Z}_2)$
and, because $H_2(F; {\Bbb Z}_2)$ and $H_2(M; {\Bbb Z}_2)$ are trivial, $\psi$ induces an isomorphism
$H_2(F ; {\Bbb Z}_2) \rightarrow H_2(\pi_1(M); {\Bbb Z}_2)$. Thus, by Theorem 5.4,
$\psi$ induces an isomorphism between the pro-2 completions of $F$ and $\pi_1(M)$.

Now $F$ admits a surjective homomorphism $\phi \colon F \rightarrow {\Bbb Z_2} \ast {\Bbb Z}_2$
sending the first free generator to the non-trivial element in the first factor,
the second free generator to the non-trivial element in the second factor,
and the remaining free generators to the identity. Composing $\phi$ with
projections onto the first and second factors gives homomorphisms $\phi_1, \phi_2 \colon F \rightarrow {\Bbb Z_2}$.
These correspond to the homomorphisms $\phi_1'$ and $\phi_2'$. Using Proposition 6.3, there is a surjective
homomorphism $\pi_1(M) \rightarrow {\Bbb Z}_2 \ast {\Bbb Z}_2$ such that
its composition with projection onto the $i$th factor is $\phi_i'$.
Theorem 3.1 gives the required surfaces $S_1$ and $S_2$, and Remark 3.2,
$[S_1] = z_1$ and $[S_2] = z_2$. $\square$

A particularly interesting case is when $M$ is the exterior of a connected finite
graph $X$ embedded in $S^3$. Then, as long as $b_1(X) > 1$, $M$ has
the same ${\Bbb Z}_2$ homology as a handlebody other
than a solid torus. We view it as quite striking that the conclusion
of Theorem 6.4 applies in this level of generality.

\vskip 18pt
\centerline{\caps 7. Further questions and remarks}
\vskip 6pt

\noindent {\caps 7.1. Making the surfaces essential}
\vskip 6pt

In 3-manifold theory, it is the surfaces that are essential that play
a particularly important role. By definition, an orientable surface
properly embedded in an orientable 3-manifold $M$ is {\sl essential}
if it is incompressible, boundary-incompressible and no component is
boundary parallel. A non-orientable surface $S$ properly embedded in $M$
is {\sl essential} if ${\rm cl}(\partial N(S) - \partial M)$ is essential.
It is well known that this has an equivalent reformulation in terms of
the way that $\pi_1(S)$ maps into $\pi_1(M)$. In particular, an essential
surface is $\pi_1$-injective.

\noindent {\bf Question 7.1.} Can one arrange for the surfaces $S_1$ and $S_2$
provided by Theorem 1.1 to be essential?

We do not have a definite answer. However, the following result establishes
that one can ensure that the surfaces are incompressible. Recall that
a (possibly non-orientable) surface $S$ properly embedded in $M$
is {\sl incompressible} if, for any embedded disc $D$ in $M$ such that
$D \cap S = \partial D$, the curve $\partial D$ bounds a disc in $S$.

\noindent {\bf Proposition 7.2.} {\sl If a compact 3-manifold $M$ contains two disjoint, properly
embedded surfaces $S_1$ and $S_2$ such that $M - (S_1 \cup S_2)$
is connected, then it contains two such surfaces that are, in addition,
incompressible.}

\noindent {\sl Proof.} We may assume that $S_1$ and $S_2$ are
both connected. Suppose that at least one is
compressible. Then, their union is compressible, by a standard
innermost curve argument. Let $D$ be a compression disc for $S_1 \cup S_2$.
Suppose that its boundary lies in $S_1$, say. Compress $S_1$ along $D$,
giving a surface $\overline S_1$. Suppose that $M - (\overline S_1 \cup S_2)$
is not connected. Now, $M - (\overline S_1 \cup S_2)$ is obtained from $M - (S_1 \cup S_2)$
by cutting along $D$ and then attaching a 2-handle. The latter operation does not
change the number of components. So, we deduce that $D$ divides $M - (S_1 \cup S_2)$
into two components, $X$ and $Y$, say. One of these components, $X$ say, lies
on the other side of $S_1$, near $\partial D$. Now, $\overline S_1$ must have two
components. This is because, near $D$, one of the parts of $\overline S_1$
has $X$ on both sides, whereas the other has $X$ on one side and $Y$ on the other.
Discard the latter component of $\overline S_1$, and let $S'_1$ be the resulting surface.
Then $M - (S_1' \cup S_2)$ is connected.
Hence, continuing in this fashion, we end up with two properly embedded, disjoint,
incompressible surfaces $S''_1$ and $S''_2$ such that $M - (S''_1 \cup S''_2)$ is
connected. $\square$

The difficulty in answering Question 7.1 in the affirmative is that
if a non-orientable properly embedded surface fails to be essential,
then there is no obvious modification that one can make to it which,
in a suitable sense, simplifies it.

\vskip 6pt
\noindent {\caps 7.2. More than two surfaces}
\vskip 6pt

This paper has been devoted to the study of two disjoint surfaces properly
embedded in a 3-manifold. It is natural to ask the following:

\noindent {\bf Question 7.3.}
Under what circumstances does a compact
orientable 3-manifold $M$ contain disjoint properly embedded surfaces 
$S_1, \dots, S_n$ such that $M - (S_1 \cup \dots \cup S_n)$ is connected, for
$n \geq 3$?

The methods in this paper do not obviously apply when $n \geq 3$.
When $n = 2$, the group ${\Bbb Z}_2 \ast {\Bbb Z}_2$ plays the central role.
This group is virtually abelian, and it is essentially for this reason that
the existence of a
surjection from a finitely generated group $G$ to ${\Bbb Z}_2 \ast {\Bbb Z}_2$ can be detected
by the pro-2 completion $\hat G_{(2)}$. However, for $n \geq 3$,
$\ast^n {\Bbb Z}_2$ is virtually free non-abelian, and so it seems
unlikely that one can detect whether a group $G$ surjects onto 
$\ast^n {\Bbb Z}_2$ purely by examining $\hat G_{(2)}$. 
This probably implies that Question 7.3 has no straightforward answer.

\vskip 18pt
\centerline{\caps References}
\vskip 6pt

\item{1.} {\caps S. Harvey}, {\sl On the cut number of a 3-manifold}, 
Geom. Topol. 6 (2002) 409--424

\item{2.} {\caps J. Hillman}, {\sl Spanning links by nonorientable surfaces.} 
Quart. J. Math. Oxford Ser. (2) 31 (1980), no. 122, 169--179.

\vfill\eject
\item{3.} {\caps M. Lackenby,} {\sl Large groups, Property $(\tau)$ and the homology growth of subgroups,}
Math. Proc. Camb. Phil. Soc. 146 (2009) 625--648.

\item{4.} {\caps M. Lackenby}, {\sl Detecting large groups}, J. Algebra 324 (2010) 2636--2657.

\item{5.} {\caps J. Stallings}, {\sl Homology and central series of groups}, 
J. Algebra 2 (1965) 170--181.

\vskip 12pt
\+ Mathematical Institute, University of Oxford, 24-29 St Giles', Oxford OX1 3LB, United Kingdom. \cr

\end